\documentclass[11pt]{amsart}
\usepackage[linktocpage]{hyperref}
\usepackage{amssymb, paralist, xspace, graphicx, url, amscd, euscript, mathrsfs,stmaryrd,epic,eepic,color}
\usepackage[all]{xy}
\SelectTips{cm}{}

\numberwithin{equation}{section}

1

\numberwithin{subsection}{section}

\allowdisplaybreaks[1]

\newtheorem*{namedtheorem}{\theoremname}
\newcommand{\theoremname}{testing}

\newtheorem{maintheorem}{Theorem}

\newtheorem{theorem}[subsection]{Theorem}

\newtheorem{proposition-definition}[subsection]
{Proposition-Definition}
\newtheorem{corollary}[subsection]{Corollary}
\newtheorem{lemma}[subsection]{Lemma}

\theoremstyle{definition}

\newtheorem{remark}[subsection]{Remark}

\newtheorem*{pf}{Proof} 

\theoremstyle{remark}

\newcommand\cA{\mathcal{A}}

\newcommand\cG{\mathcal{G}}

\newcommand\cJ{\mathcal{J}}

\newcommand\cL{\mathcal{L}}

\newcommand\cN{\mathcal{N}}
\newcommand\cO{\mathcal{O}}

\newcommand\cT{\mathcal{T}}

\newcommand\cW{\mathcal{W}}
\newcommand\cX{\mathcal{X}}
\newcommand\cY{\mathcal{Y}}
\newcommand\cZ{\mathcal{Z}}

\newcommand\CC{\mathbb{C}}

\newcommand\PP{\mathbb{P}}
\newcommand\QQ{\mathbb{Q}}

\newcommand\ZZ{\mathbb{Z}}

\newcommand\fL{\mathfrak{L}}

\newcommand\fX{\mathfrak{X}}

\newcommand\arr{\ifinner\to\else\longrightarrow\fi}

\def\displaytimes_#1{\mathrel{\mathop{\times}\limits_{#1}}}

\def\displayotimes_#1{\mathrel{\mathop{\bigotimes}\limits_{#1}}}

\newdir{ >}{{}*!/-5pt/@{>}}

\newcommand\doublelong[2]{\mathbin{\xymatrix{{}\ar@<3pt>[r]^{#1}
\ar@<-3pt>[r]_{#2}&}}}

\newlength{\ignora}

\newcommand{\e}{\acute{\hbox{e}}}

\theoremstyle{plain}
\newtheorem{thm}[subsection]{Theorem}

\newtheorem{lem}[subsection]{Lemma}

\theoremstyle{definition}

\newtheorem{rem}[subsection]{Remark}

\begin{document}
\title{Sections of Calabi-Yau threefolds with K3 fibration}

\author{Zhiyuan Li}

\address{Department of Mathematics\\
Rice University\\
6100 Main Street\\
Houston, TX 77005\\
U.S.A.}
\email{zhiyuan.li@rice.edu}


\begin{abstract}
We study sections of a Calabi-Yau threefold fibered over a curve by K3 surfaces. We show that there exist infinitely many isolated sections on certain K3 fibered Calabi-Yau threefolds and the subgroup of the N$\e$ron-Severi group generated by these sections is not finitely generated. This also gives examples of $K3$ surfaces over the function field $F$ of a complex curve with Zariski dense $F$-rational points, whose geometric model is Calabi-Yau. \end{abstract} 
\maketitle

\begin{section}{Introduction}

Let $Y$ be a smooth projective variety over the function field $F$ of a smooth projective curve $B$.  Let $Y(F)$ denote the set of $F$-rational points (or equivalently, sections of a projective  model $\cY\rightarrow B$ with generic fiber $Y$).   

The Zariski density  or potential density (i.e. density after a finite field extension) of $Y(F)$  is expected to relate to the global geometry of $Y$.  For instance,  it is known \cite{GHS} that $Y(F)$ is Zariski dense for rationally connected varieties $Y$. On the other hand, Lang's conjecture over function fields \cite{Hi},  confirmed in the curve case, predicts  that potential density fails  for certain classes of general type varieties.

For the intermediate case where the canonical class $K_Y$ is trivial,  it is generally expected that $Y(F)$ is  potentially dense. This is  known for  abelian varieties  and  $K3$ surfaces  with additional structures  \cite{BT}.  Already,  the  case for general $K3$ surfaces remains widely open.

If $F=\CC(t)$,  Hassett and Tschinkel \cite{HT} proved that a pencil of  low degree K3 surfaces  has a Zariski dense collection of sections via the study of the  deformation of sections in the pencil. 
The examples in \cite{HT} correspond to lines in the  Hilbert scheme of  the  K3 surface. It is natural to consider rational curves in that Hilbert scheme of  higher degree.  

In this paper,  we will consider  families of K3 surfaces which correspond to  conics in the Hilbert scheme. One concrete example is  the  bidegree $(2,4)$ hypersurface $X\subset \PP^1\times \PP^3$, which is a family of quartic surfaces over a conic.  Note that $X$ is a Calabi-Yau threefold,  so the sections in $X$  are generally expected to be {\it isolated } or {\it infinitesimally rigid}, i.e. the space of embedded deformations is reduced  and zero dimensional. 
Our main result is the following:
\begin{maintheorem}\label{mthm1}
There exist   countably many  isolated sections $\ell_n$  on a general  bidegree $(2,4)$ hypersurface $X\subset \PP^1\times \PP^3$ with respect to the projection to $\PP^1$.  Furthermore, the union of these sections  is Zariski dense in $X$. \end{maintheorem}

Let $NS^2(X)$ denote the N$\e$ron-Severi group  of  algebraic cycles of codimension two modulo algebraic equivalence.  It is known that this group is not finitely generated for general Calabi-Yau threefolds \cite{Vo2}.  An interesting question is to study the subgroup of $NS^2(X)$ generated by these sections.  Here we follow the  method of Clemens \cite{C3} to obtain the following theorem, which will imply the Zariski density of $\{\ell_n\}$.  \begin{maintheorem}\label{mthm2}
The subgroup $\cA$ of $NS^2(X)$ generated by  the sections is not finitely generated. 
\end{maintheorem}

This paper is organized  as follows: In section 2, we recall the N$\e$ron model theory on degenerations of intermediate Jacobians.  In particular, we state a theorem describing the N$\e$ron models coming from geometry.  Section 3 is devoted to showing the existence of  infinitely many isolated sections on $X$. We will describe the construction of these section using specialization. 
In section 4, we find a useful  degeneration of our Calabi-Yau threefolds and study the  the deformation theory of curves on the singular fiber of the degeneration.  As an application of  the result in $\S$2, we compute the group of components of the N$\e$ron model associated to the degeneration.  The main theorems are proved in section 5 and section 6. In section 7, we extend our results to higher dimensional cases.

{\bf Acknowledgements}. The  author was supported by NSF grant 0901645. The author would like to thank his advisor Brendan Hassett for introducing him to this problem, and giving lots of encouragement. Thanks Izzet Coskun for very helpful suggestions on this work, and Mattew Kerr, Zhiyu Tian, and Brian Lehmann  for helping reviewing the paper.  \end{section}

\begin{section}{Preliminaries on N$\e$ron models}
In this section, we briefly  review  some results \cite{GGK}  of N$\e$ron model theory for families of intermediate Jacobians coming from a variation of Hodge structure (VHS),   which will be used later in this paper.  For simplicity,  our VHS  arises from geometry and is paramatrized by a complex disc. 


\subsection{Geometric setting.}   Let $X$ be a smooth projective variety of dimension $2k-1$.   The intermediate Jacobian  $J(X)$ of $X$ is a compact torus defined as 
$$J(X)=H^{2k-1}(X,\CC)/(F^kH^{2k-1}(X)\oplus H^{2k-1}(X,\ZZ))$$ 
where $F^\cdot H^{2k-1}(X)$ is the Hodge filtration of $H^{2k-1}(X)$. 

More generally,   let $\Delta$ be a complex disc and let  $\pi:\cX\rightarrow \Delta$ be a  {\it  semistable degeneration},  that is:  \begin{description}
  \item[(1)] $\cX$ is smooth of dimension $2k$;
  \item[(2)] $\pi$ is projective, with the restriction $\pi:\cX^\ast=\cX\backslash\pi^{-1}(0) \rightarrow \Delta^{\ast}$ smooth, where $\Delta^\ast=\Delta\backslash\{0\}$; 
  \item[(3)] the fiber $ \cX_0=\pi^{-1}(0)$ is reduced with non-singular components crossing normally;  write $\cX_0=\cup X_i$.
\end{description} 
Consider the  VHS associated to the $(2k-1)$th cohomology along the fibres of  $\pi:\cX\rightarrow \Delta^\ast$;  then there is family of intermediate Jacobians  
 \begin{equation}\label{eq01} \cJ\rightarrow \Delta^\ast, \end{equation}
which forms  an analytic fiber space  with fiber $\cJ_s=J(\cX_s), s\in \Delta^\ast$. 

Because of the semistability assumption,  the  Monodromy theorem \cite{La} implies that the monodromy transformation $$T: H^{2k-1}(\cX_s,\ZZ)\rightarrow H^{2k-1}(\cX_s,\ZZ)$$ is unipotent.  In this situation,   Green, Griffiths and Kerr \cite{GGK} have constructed a slit analytic space  $\bar{\cJ}(\cX)\rightarrow \Delta$  such  that
 \begin{itemize}
 \item the restriction $\bar{\cJ}(\cX)|_{\Delta^\ast}$ is $\cJ\rightarrow \Delta^\ast$; 
 \item  every {\it admissible normal function} (ANF)  extends to a holomorphic section of  $\bar{\cJ}(\cX)\rightarrow \Delta$; here an ANF is a holomorphic section of (\ref{eq01}) satisfying the admissibility condition (cf. \cite{Sa} or \cite{GGK} II.B).
  \item the fiber $\bar{\cJ}_0(\cX)$ inserted over the origin   fits into an exact sequence \begin{equation}\label{eqG}0\rightarrow  \cJ_0\rightarrow \bar{\cJ}_0(\cX)\rightarrow G\rightarrow 0,\end{equation} where $G$ is  a  finite abelian group and  $\cJ_0$ is a connected, complex Lie group, considered as  the identity component of $\bar{\cJ}_0(\cX)$.
 \end{itemize}
The total space $\bar{\cJ}(\cX)$ is called the N$\e$ron model associated to $\cX$.

 \begin{remark}
In fact,  every ANF without {\it singularities} \cite{GG} extends to the identity component (cf.  \cite{GGK} II. A).\end{remark}

\begin{remark}
Kato, Nakayama and Usui have  an alternate approach  constructing N$\e$ron models  via a log mixed Hodge theory\cite{KNU}, which is homeomorphic to the construction in \cite{GGK}. (cf. \cite{Ha})
 \end{remark}


\subsection{Abel-Jacobi map}Let $\hbox{CH}^k(X)_{hom}$  be the subgroup of the Chow group of $X$ consisting of codimension $k$ algebraic cycles which are  homologically equivalent to zero.  There is an  Abel-Jacobi map   \begin{equation}
\hbox{AJ}_{X}: CH^k(X)_{hom}\rightarrow J(X)
\end{equation}
introduced by Griffiths \cite{G1}. 

Returning to the semistable degeneration $\cX\rightarrow\Delta$,   given a codimension $k$ algebraic cycle $\cZ\subset \cX$  with $Z_s=\cZ\cdot \cX_s \in CH^k(\cX_s)_{hom}$ for  $s\neq 0$,  there is an associated  admissible normal functions $\nu_\cZ$ via the Abel-Jacobi map
 \begin{equation}\label{eq02}\nu_\cZ(s)=\hbox{AJ}_{\cX_s}(Z_s ), ~s\in \Delta^\ast. \end{equation}  (cf. \cite{GGK} III)

Furthermore,  the associated function $\nu_\cZ$ will extend to the identity component of $\bar{\cJ}(\cX)$ if $\cZ$ is cohomological to zero in $\cX$.


\subsection{Threefold case.} 
With the notation above, now we assume that $\cX\rightarrow \Delta$ is a  semistable degeneration of  projective threefolds, and  denote by  $\bar{\cX}$ a  smooth projective  variety containing  $\cX$ as an analytic open subset.  

In this situation,   we have a  precise description of the group of components $G$ via  an intersection computation. \begin{theorem}\label{thm1}{(\cite{GGK} Thm.III. C.6)}
For  any multi-index $I=(i_0,\ldots, i_m)$, $|I|=m+1$, let  \begin{equation} \begin{aligned}Y_I=&\bigcap\limits_{i\in I} X_i \\  Y^{[m]}=& \coprod\limits_{|I|=m+1}Y_I . \end{aligned}\end{equation} 
Assuming that all the cohomology groups of $Y^{[m]}$ are torsion free,  then the natural map $j:Y^{[0]}\rightarrow \bar{\cX}$ induces  a sequence of  maps $$ H_4(Y^{[0]}, M)\xrightarrow{j_\ast^M} H_4(\bar{\cX}, M)\cong H^4(\bar{\cX},M)\xrightarrow{j^\ast_M} H^4(Y^{[0]}, M )\cong H_2(Y^{[0]}, M)$$  where $M=\ZZ$ or $\QQ$, and 
the composition gives  the morphism  $$\mu_M: \bigoplus\limits_{i=1}^m H_4(X_i,M) \rightarrow \bigoplus\limits_{i=1}^m H_2(X_i, M).$$
Then there is an identification of the group $G$ in (\ref{eqG}), \begin{equation}\label{eq03}G=\frac{(\hbox{Im} ~\mu_\QQ)_\ZZ} {\hbox{Im} ~\mu_\ZZ} .\end{equation} Furthermore,  the extension of the admissible normal function $\nu_{\cZ}$ (\ref{eq02}) maps to the component  corresponding to the  class $[Z_0]$ in $\bigoplus\limits_{i=1}^{m}H_2(X_i,\ZZ)$.  
\end{theorem}

\begin{remark}
 A similar result holds for a degeneration of curves. But when $\dim \cX_s>3$,  the identification $(\ref{eq03})$ may fail (cf. \cite{GGK}). 
\end{remark}
~
\end{section}

\begin{section}{Construction of  sections on K3-fibered  Calabi-Yau threefolds}

In this section, our aim is to show the existence of isolated sections on  a general  bidegree $(2,4)$ hypersurface in $\PP^1\times \PP^3$
with respect to the projection to $\PP^1$.
We begin with the construction of a hypersurface $X_0$  with at worst  nodes as singularities containing infinitely many isolated sections.  

\begin{lem}\label{lem1}
There exists a hypersurface  $X_0\subset \PP^1\times \PP^3$  with finitely many nodes, such that $X_0$ admits an infinite collection of  sections $\{\ell_{n}\} $ with respect  to the projection $X_0\rightarrow \PP^1$.  Moreover, each  $\ell_n$ lies in the smooth locus of $X_0$ and is  infinitesimally rigid. 
\end{lem}
 \begin{pf}
Let  $S\rightarrow \PP^1$ be  a  smooth rational elliptic surface, obtained by blowing up  $\PP^2$ along  nine base points  of  a general pencil of cubic curves.   This was first studied by Nagata \cite{Na},   who showed  there are infinitely many exceptional curves of the first kind,   and each of them yields  a section  $\ell_n$ of $S\rightarrow \PP^1$. 

We have a natural embedding $S\hookrightarrow  \PP^1\times \PP^2$ and  choose a smooth surface  $H \subset \PP^1\times \PP^2$ of bidegree $(1,1)$ meeting $S$ transversally in $\PP^1\times \PP^2$.   

 Let $x=(t_0,t_1; x_0,\ldots,x_3)$ be the coordinates of
$\PP^1\times \PP^3$. Consider $\PP^1\times \PP^2$ as a hyperplane of $\PP^1\times \PP^3$ defined by $x_3=0$.  Let  $| \fL| $ be the linear system of bidegree $(2,4)$ hypersurfaces in $\PP^1\times \PP^3$ containing $S$ and $H$.  Then a general member in $|\fL|$ will be a singular hypersurface with finitely many nodes contained in $S\cap H$.

More explicitly, assume that  $S$  is defined by $q(x)=x_3=0 $, while $H$ is given by the equations $l(x)=x_3=0$  for some  polynomial $q(x)$ of bidegree $(1,3)$ and $l(x)$ of bidegree $(1,1)$.    

Then a hypersurface   
 $X_0\in |\fL|$ is   given by an equation
\begin{equation} l(x)q(x)+x_3 f(x)=0\end{equation} for some bidegree
$(2,3)$ polynomial $f(x)$. The singularities of $X_0$
are eighteen nodes defined by
\begin{equation}l(x)=f(x)=x_3=q(x)=0.\end{equation}
for a generic choice of $f(x)$ by Bertini's theorem. 

For each $n$, the space of $X_0$ containing a node on $ \ell_n$  is only  a finite union of hypersurfaces in $|\fL|$.  Then we can ensure no node of $X_0$ lies on $\{\ell_n\}$ for a generic choice of $f(x)$ avoiding countably many hypesurfaces in $|\fL|$.

Furthermore, since $\ell_n^2=-1$ in $S$,  then $\cN_{\ell_n/S}=\cO_{\ell_n}(-1)$. Then the normal bundle exact
sequence \begin{equation}\label{eq21} 0\rightarrow \cO_{\ell_n}(-1)\rightarrow \cN_{\ell_n/X_0}
\rightarrow \cO_{\ell_n}(-1)\rightarrow 0\end{equation} implies that  $$\cN_{\ell_n/X_0}\cong \cO_{\ell_n}(-1)\oplus \cO_{\ell_n}(-1).$$  This proves  the infinitesimally rigidity.  \qed
\end{pf}

The following result  follows from the above lemma  and deformation theory. 
\begin{thm}\label{thm2}
For a general bidegree $(2,4)$ hypersurface  $X$ in $\PP^1\times \PP^3$, there exist  infinitely many sections $\{\ell_n\}$ on $X$ with respect to the projection $X\rightarrow \PP^1$ such that $\ell_n$ is infinitesimally rigid in $X$.
\end{thm}
\begin{pf}  From the above lemma, the rational curves $\ell_n$ in $X_0$ are stable under deformations by the  Kodaira stability theorem \cite{Ko}. This implies that the relative Hilbert scheme parameterizing the pair $(\ell,X), \ell\subset X$ is smooth over the deformation space of $X$ at $(\ell_n, X_0)$, and hence dominating.  These sections $\ell_n$ deform to nearby neighborhoods of $X_0$.  Although  $X_0$ is singular, but  we can restrict everything to the smooth locus of $X_0$ to ensure the  argument still applies. 

Furthermore, the fibration $\pi: X_0\rightarrow \PP^1$ is given by the linear system $|\pi^\ast \cO_{\PP^1}(1)|$.   Since $\pi^\ast \cO_{\PP^1}(1)$ has no higher cohomology,  it  deforms with  $X_0$ and the dimension of $|\pi^\ast \cO_{\PP^1}(1)|$ is constant by semicontinuity.  Thus the fibration will be preserved under deformation.  Note that the deformation of $\ell_n$ meets the generic fiber of  $\pi$  at one point; it follows that the deformation of $\ell_n$ remains to be  a section in a general deformation of $X_0$.  \qed
\end{pf} 

 Throughout this paper, by abuse of the notation, we continue to denote $\ell_n \subset X$  by the  section  obtained from the deformation of $\ell_n\subset X_0$.

\begin{rem}
Our method constructs infinitely many isolated rational curves of bidegree $(1,d)$ on a K3-fibered Calabi-Yau threefold in $\PP^1\times \PP^N$. For examples,  the  exceptional divisors of $S$ give degree $0$ sections on $X_0\rightarrow \PP^1$, which are of type $(1,0)$. See \cite{EJS} for  the existence of  isolated rational curves  of bidegree $(0,d)$ on K3-fibered Calabi-Yau threefolds in $\PP^1\times \PP^N$ for every integer $d\geq 1$. 
\end{rem}
~

\end{section}


\begin{section}{The degeneration of Calabi-Yau threefolds }
 In this section, we will study the degeneration of our Calabi-Yau threefolds and the deformation theory of  sections on the degenerations. 
\subsection{ An important degeneration.} 
\begin{lemma}\label{lem2}
Let $X$ be a generic bidegree $(2,4)$ hypersurface of $\PP^1\times \PP^3$.  Then there exists a projective family of bidegree $(2,4)$ hypersurfaces $\fX\rightarrow B$, containing  $X$ as a generic fiber,   such that 
\begin{itemize}
\item $\fX$ is smooth, and the generic fiber of  $\fX\rightarrow B$ is smooth;
  \item $\forall n_0\in\ZZ$, there exists a point $b_{n_0}\in B$  such that the fiber $X_{n_0}:=\fX_{b_{n_0}}$ is singular with only finitely many nodes  and satisfies \begin{enumerate}[(a)]
  \item the specialization $\ell_{n_0}\subset X_{n_0}$
passes through exactly one node  while  other specializations $\ell_n\subset X_{n_0}$ do not pass through any nodes
 for $n\neq n_0$; 
  \item  all  $\ell_{n}\subset X_{n_0}$ are  infinitesimally rigid.
\end{enumerate} 
\end{itemize}

(The notation $\fX$  is different from  $\cX$ in $\S$2.)
\end{lemma}
\begin{pf}  Consider $X$ as a deformation of the $X_0$ constructed in Lemma \ref{lem1}, where $\ell_n$ does not meet  singular locus of $X_0$.  Let $|\fL'|$ be the linear system of bidegree $(2,4)$ hypersurfaces containing $S$. The idea of the proof comes from an observation  that the space of bidegree $(2,4)$ hypersurfaces satisfying condition $(a)$ is a divisor in $|\fL'|$. 

Indeed, we can give an explicit construction as in  \cite{C1}.  With the notation in Lemma \ref{lem1},  we first consider the one parameter family of bidegree (2,4) hypersurfaces defined by  the equation \begin{equation}
l_u(x)q(x)+x_3f(x)=0,
\end{equation}
where  $
l_u(x)=u_0l_0(x)+u_1 l_1(x), u\in \PP^1$ defines  a  linear pencil of bidegree $(1,1)$ hypersurfaces. 

Let $C$ be the curve defined by \begin{equation}\label{eq30}q(x)=f(x)=0, \end{equation} meeting $\ell_n$ transversally at distinct points for a generic choice of $f(x)$. Then one can choose $l_u(x)$ outside  a countable union of
hypersurfaces in the space of all pencils, such that the hyperplane $l_u(x)=0 $ meets $S$ transversely and  does not contain more than one point of the countable set $$C\cap \left(\bigcup\limits_n \ell_n\right).$$

Then the  two parameter family  \begin{equation}\label{eq31b}\fX=\{ l_u(x)q(x)+x_3f(x)+\lambda F(x)=0 \} \rightarrow   \PP^1\times  \Delta \end{equation} 
will be the desired degeneration  for generic $F(x)$. According to our construction, for each integer $n_0$, one can find a point $u_{n_0}\in\PP^1$ such that $\fX_{(u_{n_0}, 0)}$ satisfies  condition $(a)$.

To complete the proof, it remains to show that all $\ell_n$ are infinitesimal rigid in $X_{n_0}$. When $n\neq n_0$, the rigidity of $\ell_n$ comes from  the same argument in the proof of Lemma \ref{lem1}.

If  $n=n_0$, let $X'_{n_0}$ be the blow up of $X_{n_0}$ along $P$,  and $X''_{n_0}$ the blow up of $X_{n_0}$ along $S$. Note  that  $X'_{n_0}$ and $X''_{n_0}$ are  the two  small resolutions of $X_{n_0}$. It suffices to show that  the strict transforms $\ell'_{n_0}$ and $\ell''_{n_0}$ of $\ell_{n_0}$ in $X'_{n_0}$  and $X''_{n_0}$, respectively, are infinitesimally rigid. 

Note that $\ell'_{n_0}$ is still contained in $S\subset X'_{n_0}$ as an exceptional curve,  so one can conclude that \begin{equation}\label{eq31a}\cN_{\ell'_{n_0}/X'_{n_0}}=\cO_{\ell'_{n_0}}(-1)\bigoplus \cO_{\ell'_{n_0}}(-1)\end{equation}  from  the exact sequence (\ref{eq21}).

Next, if one can find a special case of $X_{n_0}$ such that 
\begin{equation}\label{eq31}\cN_{\ell''_{n_0}/X''_{n_0}}=\cO_{\ell''_{n_0}}(-1)\bigoplus
\cO_{\ell''_{n_0}}(-1), \end{equation} then  semicontinuity will ensure that (\ref{eq31}) holds for the generic case. The existence of  such $X_{n_0}$ is  known by  Lemma 9 in \cite{Pa},   which completes the proof. \qed

\end{pf} 

\subsection{ Deforming the section through a node}
With the notation from the previous section,  let $\pi: \cX\rightarrow \Delta$ be the restriction $\fX|_{\{u_{n_0}\}\times \Delta}$, whose central fiber is  $\pi^{-1}(0)=\cX_0=X_{n_0}$.

If $m\neq n_0$,   we  know that the section $\ell_m\subset X_{n_{0}}$  deforms to a section $\ell_m(s)$ of $\cX_s$  and hence yields  a codimension two cycle $\cL_m\subset \cX$ with   \begin{equation}\label{eq32}\cL_m\cdot \cX_s=\ell_m(s),~s\in \Delta. \end{equation}

However,  $\ell_{n_0}$  in $X_{n_0}$  cannot deform with $X_{n_0}$  in $\cX$ since  there is a non-trivial obstruction for first order deformations.   This obstruction will vanish after  a degree two base change. In this subsection, we will show that this is a sufficient condition to deform $\ell_{n_0}$ with $X_{n_0}$.   The following result  is inspired by  \cite{C3}.  
 
\begin{theorem}\label{thm3}
The section $\ell_{n_0}\subset X_{n_0}$  can deform with $X_{n_0}$ in  $\cX$ only  after a degree two base change. In other words, we have the following diagram

\begin{equation}\label{eqthm}
\xymatrix{  \mathcal{L}_{n_0} \ar[d] \ar@{^{(}->}[r]^i &\widetilde{\cX} \ar[d]_{\tilde{\pi}} \ar[r]& \cX \ar[d]^{\pi} \\ 
\widetilde{\Delta}\ar[r]^{Id}& \widetilde{\Delta} \ar[r]^d & \Delta}
\end{equation}
where $d:\widetilde{\Delta}\rightarrow \Delta$ is the double
covering map of the disc $\Delta $ ramified at the center $0\in
\Delta$, and $\mathcal{L}_{n_0}\cap \tilde{\pi}^{-1}(0)=\ell_{n_0}$.
\end{theorem}
Before proceeding to the proof, let us fix some notation as follows:
\begin{itemize}
  \item $X_{n_0}$ is defined by the equation  $F_0(x)=0$, without loss of generality, having  a node   $p_0=(1,0; 1,0,0,0)\in \PP^1\times \PP^3$ ; 
  \item the section $\ell_{n_0}\subset X_{n_0}$ passing through  $p_0$ is parametrized by the morphism 
  \begin{equation}
 \begin{aligned}
    \phi: \PP^1&\longrightarrow X_{n_0}~~~~~~~~~~~~~~~~~~
\\
    t=(t_0,t_1) & \mapsto (t_0,t_1;  \phi_0(t),\ldots, \phi_3(t))
 \end{aligned}
\end{equation}
with $\phi(1,0)=p_0$ for some  degree $d$ homogenous polynomials $\phi_i(t)$, \\$ i=0, \ldots, 3$;
  \item  the family $\cX\rightarrow \Delta$  is given by the equation  \begin{equation}
F_0(x)+  sF(x)=0,~ s\in \Delta,
\end{equation} 
for some polynomial $F(x)$, with $F(p_0)\neq 0$. 
\end{itemize}  

\noindent With the notation above,  we  first give an explicit description of the global sections of the normal sheaf $\cN_{\ell_{n_0}/X_{n_0}}$.  
\begin{lemma}\label{lem3}
A global section of $\cN_{\ell_{n_0}/X_{n_0}}$ can be represented  by a set of homogenous polynomials
\begin{equation}\label{eq32} \{(\sigma_i(t))_{i=0,1},
(\delta_j(t))_{j=0,\ldots,3}\}\end{equation} with $deg(\sigma_i)=1$ and $deg(\delta_j)=d$, 
subject to  the condition
\begin{equation} \sum\limits_{i=0}^1\sigma_i(t)\frac{\partial F_0}{\partial
t_i}(\phi(t)) +\sum\limits_{j=0}^3 \delta_j(t)\frac{\partial
F_0}{\partial x_j}(\phi(t))=0.
\end{equation}

Moreover, (\ref{eq32}) is a trivial section of $\cN_{\ell_{n_0}/X_{n_0}}$ if and only if it 
satisfies the condition
\begin{equation}\label{eq33} \delta_j(t)=\sigma_0(t)\frac{\partial
\phi_j}{\partial t_0}+\sigma_1(t)\frac{\partial\phi_j}{\partial
t_1},~j=0,\ldots, 3.\end{equation}

\end{lemma}

\begin{pf} Let us denote the invertible sheaf $\pi^\ast_1\cO_{\PP^1}(a)\otimes \pi^\ast_2\cO_{\PP^3}(b)$ by $\cO_{\PP^1\times\PP^3}(a,b)$, where $\pi_1$ and $\pi_2$ are natural projections of $\PP^1\times \PP^3$. 
Let $\cT_X$ be the tangent sheaf of $X$. Due to the exact sequence
\begin{equation}
0\rightarrow \cO_{X_{n_0}}^{\oplus 2} \rightarrow \cO^{\oplus 2}_{X_{n_0}}(1,0)\oplus \cO^{\oplus 4}_{X_{n_0}} (0,1) \rightarrow \cT_{\PP^1\times \PP^3}|_{X_{n_0}}\rightarrow 0
\end{equation}
and 
\begin{equation}
0\rightarrow \cT_{X_{n_0}}\rightarrow \cT_{\PP^1\times \PP^3}|_{X_{n_0}}\rightarrow \cO_{X_{n_0}}(2,4)\rightarrow 0
\end{equation} 
one can express  a global section of $\cT_{X_{n_0}}$ as a set of bidegree homogenous polynomials $\{(\sigma_i)_{i=0,1}; (\delta_j)_{j=0,\ldots,3}\}$ satisfying 
\begin{equation}
 \sum\limits_{i=0}^1\sigma_i\frac{\partial F_0}{\partial
t_i}+\sum\limits_{j=0}^3 \delta_j\frac{\partial
F_0}{\partial x_j}=0,
\end{equation}
 where $\sigma_i$ are of  bidegree $(1,0)$,  while $\delta_j$ are of bidegree $(0,1)$. 
 
Then the statement  follows from the following exact sequence, 
\begin{equation}\cT_{\ell_{n_0}}\rightarrow\cT_{X_{n_0}}|_{\ell_{n_0}}\rightarrow \cN_{\ell_{n_0}/X_{n_0}}\longrightarrow
0 \end{equation}
where the induced map  $g:H^0(\ell_{n_0}, \cT_{\ell_{n_0}})\rightarrow H^0(\ell_{n_0}, \cT_{X_{n_0}}|_{\ell_{n_0}})$ can be expressed as 
\begin{equation}a_0\frac{\partial }{\partial t_0}+a_1\frac{\partial }{\partial
t_1}\longmapsto (a_0,a_1; a_0\frac{\partial \phi_i}{\partial
t_0}+a_1\frac{\partial \phi_i}{\partial t_1})_{i=0,\ldots,
3}\end{equation} \qed
\end{pf}

\noindent{\bf Proof of Theorem \ref{thm3}}.  Let  us make the base change  $\tilde{\Delta}\rightarrow\Delta$ sending $r$ to $r^2$, and write
  \begin{equation}\label{eq34}\tilde{\cX}:=\{  F_0(x)+r^2 F(x)=0, ~~r\in \tilde{\Delta}\}. \end{equation}
To prove the assertion,  it suffices to show the existence of  a formal deformation $\Phi(r,t)$ in (\ref{eq34}),  i.e.  there is  a  sequence of maps \begin{equation}\phi^{[k]}(t)=(t; \phi_0^{[k]}(t), \ldots, \phi_3^{[k]}(t) )\in\PP^1\times \PP^3, ~k\geq 0,\end{equation} with $deg(\phi^{[k]}_i(t))=d $ and  $\phi^{[0]}=\phi$,  such that the power series  \begin{equation} \Phi(r,
t)=(t_0,t_1;  \sum\limits_{k=0}^{\infty} r^k \phi_i^{[k]}(t))_{i=0,1\ldots,3} \end{equation} satisfies the condition
\begin{equation}\label{eq35}
F_0(\Phi(r,t))+r^2 F(\Phi(r,t))=0.\end{equation}
Our proof of the existence of  $\Phi(r,t)$ proceeds as follows:
~

\noindent{(I) \it First order deformation.}  The first order deformation of $\phi$ is determined by $\phi^{[1]}(t)$, which can  be solved   by differentiating
(\ref{eq35}) with respect to $r$  and setting $r=0$. Hence we obtain
\begin{equation}\label{eq36}
\sum\limits_{i=0}^3  \frac{\partial F_0}{\partial
x_i}(\phi(t))\phi^{[1]}_i(t)=0.
\end{equation}
Note that  (\ref{eq36}) is a homogenous  polynomial of degree $4d+2$, which has $4d+3$
coefficients and the coefficient of the $t_0^{4d+2}$ term is zero by assumption.  
Then one can consider (\ref{eq36}) as 
 $(4d+2)$ equations in $4(d+1)$ unknowns and 
 denote  $M_{(\phi, ~F_0)}$ by the $(4d+2)\times (4d+4)$ matrix corresponding to the system of these equations.  

Our first claim is that the  $M_{(\phi,~F_0)}$ is of  full rank, which  is equivalent to saying that the dimension of  the solution space of $\phi^{[1]}(t)$ is two.

By Lemma \ref{lem3},  the set   \begin{equation}\label{eq37}
\{(t_0,t_1); \phi^{[1]}_0,\ldots, \phi_i^{[1]})  \}     \end{equation}
gives  a  global section 
of $\cN_{\ell_{n_0}/X_{n_{0}}}$ and  is  trivial  if  and only if \begin{equation}\phi^{[1]}_i=t_0\frac{\partial \phi_i}{\partial
t_0}+t_1\frac{\partial \phi_i}{\partial t_1}\end{equation}
by (\ref{eq33}).  So if  rank $M_{(\phi,~F_0)} \leq 4d+1$,  then $\dim H^0(\ell_{n_0}, \cN_{\ell_{n_0} /X_{n_0}} ) \geq 2$.  

Let $ev: H^0(\ell_{n_0}, \cN_{\ell_{n_0}/X_{n_0}})\rightarrow \CC^3 $  be the evaluation map at $p_0$.  As in \cite{C3} $\S$3, one can show that there is at most one condition  lifting the analytic  section of  $\cN_{\ell_{n_0}/X_{n_0}}$ to a section of  $\cN_{\ell'_{n_0} /X'_{n_0}}$,   because  the image of $ev$ is at most two dimensional, while the composition of  the  sequence of evaluation maps at  $p_0$
 \begin{equation}\label{eq37} \cN_{\ell'_{n_0}/X'_{n_0}}\rightarrow  \cN_{\ell_{n_0}/X_{n_0}}\rightarrow \CC^3 \end{equation}
only has a one dimensional image.

However, from (\ref{eq31a}), we know 
 that  $H^0(\ell'_{n_0}, \cN_{\ell'_{n_0}/X'_{n_0}})=0$.  Thus we prove the first claim by contradiction. 
\\

\noindent{(II) \it Higher order.}  We  continue to solve  $\phi^{[2]}(t)$  by differentiating (\ref{eq35}) twice, and thus obtain 
\begin{equation}
\sum\limits_{i=0}^3  \frac{\partial F_0}{\partial
x_i}(\phi(t))\phi^{[2]}_i(t)=-\sum\limits_{i,~j} \frac{\partial F_0}{\partial x_i
\partial x_j}(\phi(t))\phi^{[1]}_i (t)\phi^{[1]}_j(t) +2F(\phi(t)) .
\end{equation}
Obviously,   there is a non-trivial  obstruction to lift  $\phi^{[1]}(t)$ to  second order given by the equation,
\begin{equation}\label{eq38}
\sum\limits_{i,~j} \frac{\partial F_0}{\partial x_i
\partial x_j}(p_0)\phi^{[1]}_i (1,0)\phi^{[1]}_j(1,0) =-2F(p_0)\neq 0.
\end{equation}
Any $\phi^{[1]}(t)$ satisfying (\ref{eq38}) can be lifted to the second order, because $M_{\phi,F_0}$ is full rank. 

Our second claim is that there exists a first order deformation which  can be lifted to  second order.  Otherwise,  every $\phi^{[1]}(t)$ will satisfy the condition
\begin{equation}\label{eq39}
\sum\limits_{i,~j} \frac{\partial F_0}{\partial x_i
\partial x_j}(p_0)\phi^{[1]}_i (1,0)\phi^{[1]}_j(1,0) = 0.
\end{equation}
From the above assumption, there is an non-trivial analytic section $\alpha \in H^0(\ell_{n_0}, \cN_{\ell_{n_0}/X_{n_0}})$, whose image via the evaluation map  at $p_0$ lies  in the tangent cone $C_{p_0}$ of  $X_{n_0}$ at $p_0$ and  is normal to the  tangent direction of $\ell_{n_0}$.  

As in \cite{C3} $\S$3, this means  that $ev(\alpha)$ lies  in the union of the images of
$\cN_{\ell'_{n_0}/X'_{n_0}}$ and  $\cN_{\ell''_{n_0}/X''_{n_0}}$ in  $\cN_{\ell_{n_0}/X_{n_0}}$. This is a contradiction, because none of the non-trivial sections of $\cN_{\ell_{n_0}/X_{n_0}}$  can lift  by Lemma \ref{lem3}.  Hence there exists $\bar{\phi}^{[1]}(t)$ satisfying (\ref{eq38}), and the second claim is  proved.   

Furthermore,   set  $b_i=\bar{\phi}_i^{[1]}(1,0)$.  It is not difficult to see that the equation (\ref{eq36}) along with 
 \begin{equation}
\sum\limits_{i,~j} \frac{\partial F_0}{\partial x_i \partial
x_j}(p_0)\phi^{[1]}_i(1,0)b_j =0
\end{equation}
only has a one dimensional set of  solutions. Hence  the associated $(4d+3)\times (4d+4) $ matrix   $M'_{\bar{\phi}^{[1]}, F_0}$ is  full rank. 

For higher orders,  $\phi^{[n]}(t)$ is determined by the equation \begin{equation}
\sum\limits_{i=0}^3  \frac{\partial F_0}{\partial
x_i}(\phi(t))\phi^{[n]}_i(t)=some~polynomial~given~by ~\phi^{[k]}~for~k<n,
\end{equation}
while the obstruction to (n+1)th order is \begin{equation}
\sum\limits_{i,~j} \frac{\partial F_0}{\partial x_i \partial
x_j}(p_0)\phi^{[n]}_i(1,0)b_j =some~number~given~by ~\phi^{[k]}~for~k<n.
\end{equation}
Then one can solve $\phi^{[n]}(t)$ by induction because of the full rank of $M'_{\bar{\phi}^{[1]}, F_0}$.  \qed

\subsection{Semistable degeneration.}
In this subsection, we will desingularize the family $\widetilde{\cX}$ to obtain a semistable degeneration, and identify the group of components associated to this semistable degeneration. 

Let  $\cW$ be the blow up of  $\tilde{\cX}$  along all the nodes on $\tilde{\cX}_0$. Then we have  
\begin{itemize}
  \item the ambient space $\cW$ is smooth, and the generic fiber  of $\cW\rightarrow \tilde{\Delta}$ is smooth;
  \item the central fiber $\cW_0=\bigcup\limits_{i=0}^{18} W_i$ is strictly normal crossing, where 
  \begin{enumerate}[(1)]
  \item $W_0$ is the blow up of $\tilde{\cX}_0$ along all the nodes;
  \item $W_i$ are disjoint smooth quadratic threefolds in $\PP^4$, meeting $W_0$ transversally at the exceptional divisor $E_i\cong \PP^1\times \PP^1$ for $i=1,\ldots, 18$.
\end{enumerate}
 \end{itemize}

As an application of Theorem \ref{thm1}, we shall apply (\ref{eq01}) to compute the group of components of the N$\e$ron model $\bar{J}(\cW)$ associated to the semistable degeneration $\cW\rightarrow \tilde{\Delta}$. 

In order to give  a geometric description of the homology groups of each component of $\cW_0$, we first set the following notation:   
\begin{itemize}
  \item $P$ is the strict transform of bidegree $(1,0)$ hyperplane section of $ X_{n_0}$ in $W_0$, and $ D $ is a generic fiber of $W_0$ over $\PP^1$; 
 $\tilde{H}$ is the strict transform of $H$ in $W_0$;  
  
 \item $E_i$ are exceptional divisors of $W_0, i=1,2,\ldots, 18$.
  \item $Q_i$ is the hyperplane section of $W_i$,  $i=1,2,\ldots, 18$.
  \item  $L$ is a line on the fiber of $W_0$ over $\PP^1$,  and $L'$ is a section of $W_0$ with respect to the projection; $C'$ is the proper transform of the curve (\ref{eq30}) in $W_0$;  $R_i$ is one of the ruling of $E_i$.
  \item $L_i$ is the line in $W_i$.
\end{itemize}
Then the integral basis for these homology groups can be represented by the fundamental class of the algebraic cycles above: 
\begin{equation}\begin{aligned} H_2(W_0)&=\langle[L], [L'], [C'],  [R_1],\ldots [R_{18}]\rangle, ~   &&H_4(W_0)=\langle[D],[K], [ \tilde{P}], [E_1],\ldots [E_{18}]\rangle;\\
H_2(W_i) &=\langle[L_i]\rangle, ~i=1,2\ldots, 18;  &&H_4(W_i)=\langle[Q_i]\rangle, ~i=1,2\ldots, 18.
\end{aligned} \end{equation}
By a straightforward computation, we can express the map $$\mu_\ZZ: \bigoplus\limits_{i=0}^{18}H_4(W_i,\ZZ)\longrightarrow \bigoplus\limits_{i=0}^{18} H_2(W_i,\ZZ)$$  as a matrix:  $$\begin{array}{c|c|c|c|c|c}~ & [P] & [D]&[ \tilde{H}] & [E_i] & [Q_j]  \\\hline[ L] & 1 & 0 & 1 & 0 & 0 \\\hline [L' ]& 0 & 1& 1 & 0 & 0 \\\hline [C'] & 0 & 0 & 18& 1 & 0\\\hline [R_k] & 0 & 0 & 1 & 2 \delta_{ik}& 2\delta_{jk} \\\hline[ L_l] & 0 & 0 & 1 &2 \delta_{il} & 2\delta_{jl}\end{array}$$
Thus the  group  is computed as \begin{equation}\label{eq310}G=\frac{\hbox{Im}(\mu_\QQ)_\ZZ}{\hbox{Im}(\mu_\ZZ)}= \frac{\bigoplus\limits_{k=1}^{18} \ZZ([R_k]+[L_k])}{\bigoplus\limits_{k=1}^{18} \ZZ (2[R_k]+2[L_k])\bigoplus \ZZ \sum\limits_{i=1}^{18} ([R_i]+[L_i]) }.\end{equation}
~

Furthermore, let $\check{\cL}_{n_0}$ denote the strict transform of $\cL_{n_0}$ in $\mathcal{W}$. The following lemma is straightforward (cf. \cite{C2}) and  will be used in next section.
\begin{lemma}\label{lem4}
Let  $E_{i_0}$ be the exceptional divisor in $W_0$ corresponding to the node which $\ell_{n_0}$ passes through, and $\check{\ell}_{n_0}$ the strict transform of $\ell_{n_0}$ in $W_0$. Then \begin{equation}\check{\cL}_{n_0}\cap W_0=\check{\ell}_{n_0}+ (\hbox{one of  the rulings of $E_{i_0}$}) .\end{equation}
\end{lemma}
~

\end{section}


\begin{section}{Infinite generation of the subgroup of Griffiths group}

\subsection{Griffiths group of Calabi-Yau threefolds} Let $X$ be a smooth projective threefold and $\hbox{CH}^2(X)_{alg} $ be the subgroup of $\hbox{CH}^2(X)$  consisting of  codimension $2$ cycles which are algebraically equivalent to zero. The Abel-Jacobi image \begin{equation}\hbox{AJ}_{X} (\hbox{CH}^2(X)_{alg})=A \subseteq J(X) \end{equation} is an abelian variety.
 The abelian variety $A$ is called the algebraic part of $J(X)$. It lies in the largest complex subtorus $J(X)_{alg}\subset J(X)$, whose tangent space at $0$ is contained in $H^{1,2}(X)$. (cf.  \cite{Vo}.vI) 

 The Griffiths group $\hbox{Griff}^2(X)$ is the quotient $\hbox{CH}^2(X)_{hom}/ \hbox{CH}^2(X)_{alg}$,  which is a subgroup of $NS^2(X)$.  
In the case of  Calabi-Yau threefolds,   the following result is mentioned in \cite{Vo1} and \cite{Vo2}.
\begin{theorem}
If $X$ is a nonrigid Calabi-Yau threefold, i.e. $h^1(\cT_{X})\neq 0$,
then $J(X_s)_{alg}=0$  for a  general deformation $X_s$ of $X$. In particular, $\mathrm{AJ}_{X_s}$ factors $$ \mathrm{AJ}_{X_s} : \mathrm{Griff}^2(X_s) \rightarrow J(X_s)$$
\end{theorem}
Now we return to the case of $X$ a generic bidegree $(2,4)$ hypersurface of $\PP^1\times \PP^3$, and hence $\hbox{AJ}_X: \hbox{Griff}^2(X)\rightarrow J(X)$ is well defined.  Recalling that the group $\cA\subset NS^2(X) $ is generated by $\{\ell_n\}$ of different degree, we consider the non-trivial group $\cA\cap \hbox{Griff}^2(X)$.

Since the rank of $H_2(X,\ZZ)$ is two by the Lefschetz hyperplane theorem, there exists integers $a, a_n,b_n$, such that \begin{equation}\psi_n:=a\ell_n-a_n\ell_0-b_n\ell_1 \equiv_{hom} 0, ~\forall n\in \ZZ.  \end{equation} 

As in Remark 3.3,  the fundamental class of $\ell_n$ is of type $(1,d_n)$ in $H_2(X,\ZZ)$ and we can assume $\ell_0$ is of type $(1,0)$. After a suitable choice of $\ell_1$,  we can select $a$ to be odd. In fact, let us denote $\kappa$ by the largest number such that $d_n$ is divisible by  $2^\kappa $ for all $n\in \ZZ$,  and denote $\ell_1$ by the section of type $(1,d_1)$ satisfying that $d_1/2^\kappa$ is  odd,  one can choose $a=d_1/2^\kappa$ as desired. 
Denote $\cG$ by the subgroup of $\cA\cap NS^2(X)$ generated by $\psi_n$. 

 \subsection{Infinite generation of $\cG$} 
 In this subsection, we shall prove the following result, which implies Theorem \ref{mthm2}. 
 
 \begin{theorem}
 The Abel-Jacobi image $\mathrm{AJ}_X(\cG)\otimes\QQ$ is of infinite rank for generic $X$. 
 \end{theorem}

\begin{pf}
Suppose that there is a relation
\begin{equation}\label{eq41}\sum\limits_{finite}  c_n \hbox{AJ}_X( \psi_n)=0 \end{equation}   for generic $X$, which gives $\hbox{AJ}_X(\sum\limits_{finite} ac_n\ell_n)=0 $.    In particular,   we  assume that (\ref{eq41}) holds for generic fiber of the two parameter family   $\fX\rightarrow \PP^1\times \Delta$ (\ref{eq31b}). 

From the construction in Lemma \ref{lem2},  for each integer $n$,  there is a point $(u_n,0)\in \PP^1\times \Delta$ such that the fiber $X_n=\fX_{(u_n,0)}$ over $(u_n, 0)$ satisfies the conditions $(a)$ and $(b)$.  Let $\cX\rightarrow \Delta$ be the restriction $\fX|_{\{u_n\}\times\Delta} $  and $\tilde{\cX}=\cX \times_{\Delta} \tilde{\Delta}$ for a degree two base change $\tilde{\Delta}\xrightarrow{r\mapsto r^2} \Delta$.  
Then  the  family of cycles \begin{equation}\cZ=ac_n\cL_n+\sum\limits_{m\neq n} ac_m\tilde{\cL}_m\end{equation}  satisfies   $\hbox{AJ}_{\tilde{\cX}_r}(\cZ\cdot \tilde{\cX}_r)=0$, where $\cL_n$ is given by (\ref{eqthm}) in Theorem \ref{thm3} and $\tilde{\cL}_m$ is the lift  of  \ref{eq32} in $\tilde{\cX}$ for $m\neq n$.

To make use of the N$\e$ron model, we blow up $\tilde{\cX}$ along all the nodes to get the semistable degeneration $\cW\rightarrow \tilde{\Delta}$ as in $\S$ 4.6. There is an associated  N$\e$ron model $\bar{J}(\cW)$. Denote by $\check{\cZ}$ the lifting of $\cZ$ in $\cW$;    then  the associated admissible normal function $\nu_{\check{\cZ}}$ is a zero holomorphic section  and naturally extends to the identity component of $\bar{J}(\cW)$.  

On the other hand,  write $\check{Z}_0:=\check{\cZ}\cdot \cW_0$;  then  $\nu_{\check{\cZ}}$ extends to the component corresponding to the class of $[\check{Z}_0]$ in  $ \bigoplus\limits_{i=0}^{18} H_2(W_i,\ZZ)$  by Theorem \ref{thm1}.   According to Lemma \ref{lem4},  we have    \begin{equation}[\check{Z}_0]=ac_n([R_n]+[L_n])+linear~ combinations ~of ~[L], [L']~ and ~[C'] \end{equation} which corresponds to $ac_n([R_n]+[L_n])$ in $G$. Then as indicated in $\S$4.6,   $\nu_{\check{\cZ}}$ extends to the identity component  if and only if  $c_n$ is even, because $a$ is odd.
 
Repeating the above process  for each integer $n$,  one proves that all the coefficients in (\ref{eq41}) are even. Thus, the elements $\{\hbox{AJ}_X(\psi_n)\}$ are linearly independent modulo two which implies that
$\hbox{AJ}_X(\cG)\otimes \ZZ_2$ has infinite rank. Then the assertion follows from the  lemma below. 
\begin{lemma}
Let $M$ be an abelian group with $$M_{torsion}\subseteq (\QQ/\ZZ)^r.
$$ Then $rank_{\ZZ_2} (M\otimes \ZZ_2 )\leq rank_\QQ (M\otimes
\QQ)+r$.
\end{lemma}
\qed
\end{pf}
\begin{remark}
Note that the monodromy of the degeneration $\cW\rightarrow \tilde{\Delta}$ satisfies $(T-I)^2=0$. One can take an alternate approach by using Clemens'  N$\e$ron model $\bar{J}^{Clemens}(\cW)$ \cite{C2} to extend the associated normal function.  Actually, the identity component of  $\bar{J}(\cW)$ is a subspace of Clemens's N$\e$ron model (cf. see also  \cite{Zu}).

\end{remark}

\begin{remark}
Let $\iota:\cW\rightarrow \cW$ be the natural involution. The admissible normal function $\nu_{\cZ'}$ associated to the family of  algebraic cycles $\cZ'=\check{\mathcal{L}_n}-\iota(\check{\mathcal{L}_n})$ extends to the identity component \cite{GGK}.  One can also prove the infinite generation of $G$   by showing that $\nu_{\cZ'}(0)$ is a nontrivial element in the identity component. The proof  is similar to our computation of the group of  components (cf. \cite{Ba}). 
\end{remark}
\end{section}


\begin{section}{Proof of the Main theorem}

\noindent{\it Proof of Theorem \ref{mthm1}}.  
Assume to the contrary that the union of the sections $\{\ell_n\}$ is not Zariski dense in $X$. Let $\Sigma$ be the Zariski closure of the union of these curves, and $\tilde{\Sigma}$  the desingularization of $\Sigma$. 
Then the  proper morphism  $$\varphi: \tilde{\Sigma}\rightarrow X$$ induces a homomorphism $$\varphi_\ast: NS^1(\tilde{\Sigma})\rightarrow NS^2(X).$$
The homomorphism $\varphi_\ast$ maps the algebraic cycle $\ell_n$ in $\tilde{\Sigma}$ to the corresponding 1-cycle in $X$. So the group $\cA$ is contained in the image of $NS^1(\tilde{\Sigma})$ via $\varphi_\ast$.

It is well known that $NS^1(\tilde{\Sigma})$ is a finitely generated abelian group by the N$\e$ron-Severi theorem, which contradicts Theorem \ref{mthm2}. This completes the proof. \qed

\begin{remark}
Our result can be generalized to other Calabi-Yau threefolds fibered by complete intersection K3 surfaces in $\PP^n$. For instance,  one can find an analogous statement for some Calabi-Yau threefolds in $\PP^1\times \PP^4$ fibered by the complete intersection of a quadratic and a cubic in $\PP^4$.   \end{remark}

\begin{remark}
Our method does not yield examples in $\bar{\QQ}(t)$, since all type $(2,4)$ hypersurfaces $X$ over $\bar{\QQ}$ might lie in a countable union of ``bad" hypersurfaces of the parameterization space.

Furthermore, one can also see \cite{Be} \cite{Bl} for conjectures of $CH_0(Y)_{hom}$ when $Y$ is a surface over a number field or a function field of a curve defined over a finite field. \end{remark}
\end{section}

\begin{section}{Higher dimensional Calabi-Yau hypersufaces in $\PP^1\times\PP^N$}
In this section,  we consider the case of bidegree $(2,N+1)$ hypersurfaces  in $\PP^1\times \PP^N$ for $N\geq 3$.  The following theorem is obtained via a similar argument as Lemma \ref{lem1}. 
\begin{theorem}
For a general hypersurface $X^N\subset \PP^1\times \PP^{N}$ of bidegree $(2,N+1)$, there exists an infinite sequence of sections $\{\ell_k\}$ on $X^{N}$ of different degrees with respect to the projection $X^N\rightarrow \PP^1$,  such that \begin{equation}\cN_{\ell_k/X^N}=\cO_{\ell_k}(-1)\oplus\cO_{\ell_k}(-1)\oplus \overbrace{\cO_{\ell_k}\oplus\ldots \oplus\cO_{\ell_k} }^{N-3}.
\end{equation}
Furthermore, the subgroup $\cG^N\subset NS^2(X^N)$ generated by the algebraic codimesion $2$-cycles $\psi_k^N$, which are swept out by the deformations of $\ell_k$, is not finitely generated.
\end{theorem}
\begin{proof} The proof will proceed by induction on $N$. Suppose our statement holds for $N=m\geq 3$. When $N=m+1$, it suffices to produce a singular one with the desired properties. The construction is as follows,

Let us denote the coordinate of $\PP^1\times \PP^{m+1}$ by $x=(s,t;x_0,\ldots, x_{m+1})$. Then we consider the bidegree $(2,m+2)$ hypersurface $X_0^{m+1}$ defined by the equation
 \begin{equation}
x_0g(x)+x_nh(x)=0
\end{equation}
for some bidegree $(2,m+1)$ polynomials $g,h$.

Now, we choose our $g(x), h(x)$ satisfying the following conditions: 
\begin{enumerate}[\bf(1)]
  \item $X^{m+1}_0$ is only  singular at \begin{equation}\label{eq81} x_0=x_{m+1}=g(x)=h(x)=0; \end{equation}
  \item the subvariety \begin{equation}X^m:=\{x_{m+1}=g(x)=0\}\end{equation} satisfies the inductive assumption. Denote by $\ell_k$ the corresponding sections on $X^m$;
  \item All the sections $\ell_k$ lie outside the singular locus (\ref{eq81}).
\end{enumerate}

Similar as in the proof of Lemma \ref{lem1},  condition (1) will be satisfied due to Bertini's theorem and condition (3) can be achieved for a generic choice of $g(x)$ outside countably many hypersurfaces of the parametrization  space of bidegree $(2,m)$ polynomials. 

Next, we compute the normal bundle $\cN_{\ell_k/X_0^{m+1}}$ from the following exact sequence:
\begin{equation}0\rightarrow \cN_{\ell_k/X^m}\rightarrow \cN_{\ell_k/X_0^{m+1}}\rightarrow \cN_{X^m/X_0^{m+1}}|_{\ell_k}\rightarrow 0.\end{equation}
By assumption, we have $\cN_{\ell_k/X^{m}}= \cO_{\ell_k}(-1)^{\oplus 2}\oplus\cO_{\ell_k}^{\oplus m-3}$. Since $$ \cN_{X^m/X_0^{m+1}}|_{\ell_k}=\cO_{\ell_k},$$ it follows that  
\begin{equation}
\cN_{\ell_k/X^{m+1}_0}= \cO_{\ell_k}(-1)^{\oplus 2}\oplus\cO_{\ell_k}^{\oplus m-2}.
\end{equation}

Moreover, let $\cG^{m+1}\subset NS^2(X^{m+1}_0)$ be the subgroup generated by deformations of $\ell_k$ in $X_0^{m+1}$.  There is a morphism  \begin{equation}i^\ast: NS^2(X_0^{m+1})\rightarrow NS^2(X^m).\end{equation} induced by the inclusion $X^m\xrightarrow{i} X_0^{m+1}$. Then the infinite generation of  $\cG^{m+1} $ follows from the inductive assumption on $X^m$. This completes the proof.
\end{proof}

As in  $\S$6,  the following result is deduced from the infinite generation of $\cG^N$:
\begin{corollary}
The sections on general bidegree $(2,N+1)$ hypersurfaces of $\PP^1\times\PP^{N}$ are Zariski dense.
\end{corollary}
\end{section}

\bibliographystyle{alpha}
\bibliography{SCTWKF}


\end{document}